\documentclass{amsart}
\usepackage{hyperref}
\usepackage{a4wide}
\usepackage{amsfonts}
\usepackage{url}
\usepackage{listings}
\usepackage[osf,sc]{mathpazo}
\theoremstyle{plain}

\usepackage{tcolorbox}

\usepackage{breqn}

\theoremstyle{definition}

\allowdisplaybreaks

\subjclass[2010]{33D99}

\keywords{Theta functions, Ramanujan's theta function, Jacobi's theta function} 

\author{Sumit Kumar Jha}
\address{G-2, C-185, Ramprastha, Ghaziabad, U.P.-201011, India}
\email{sumitkumarjha.iiit@gmail.com} 

\begin{document}
	
\title{An identity for $(-q^{5};\,q^{10})_{\infty}$}

\begin{abstract}
We prove that
$$
\prod_{n=0}^{\infty}(1+q^{10n+5}) = \frac{\sum_{n=-\infty}^{\infty}q^{n^{2}}\, \sum_{n=-\infty}^{\infty}(-1)^{n}\, (-q)^{n(3n-1)/2}}{4\, \sum_{n\geq 0}(-q)^{\frac{n(n+1)}{2}}\sin \left\{\frac{(2n+1)3\pi}{10}\right\}\, \sum_{n\geq 0}(-q)^{\frac{n(n+1)}{2}}\sin \left\{\frac{(2n+1)\pi}{10}\right\}} \qquad |q|<1
$$
using identities due to Ramanujan.
\end{abstract}

\maketitle
\section{Notations and Definitions}
We assume that $|q|<1$ hereon. The functions $\varphi(q)$, $\psi(q)$, $f(-q)$, and $\chi(q)$ are defined as
\begin{align*}
    \varphi(q)&:=\sum_{n=-\infty}^{\infty}q^{n^{2}},\\
    \psi(q)&:=\sum_{n=0}^{\infty}q^{n(n+1)/2},\\
    f(-q)&:=\sum_{n=-\infty}^{\infty}(-1)^{n}q^{n(3n-1)/2}=\prod_{n\geq 1}(1-q^{n}),\\
    \chi(q)&:=\prod_{n\geq 0}(1+q^{2n+1}).
\end{align*}
The four \emph{Jacobi theta functions} \cite[p. 464]{Watson} are defined by:
\begin{align*}
\theta_{1}\left(z,q\right)&=2\sum\limits_{%
n=0}^{\infty}(-1)^{n}q^{(n+\frac{1}{2})^{2}}\sin\left((2n+1)z\right),\\
\theta_{2}\left(z,q\right)&=2\sum\limits_{%
n=0}^{\infty}q^{(n+\frac{1}{2})^{2}}\cos\left((2n+1)z\right),\\
\theta_{3}\left(z,q\right)&=1+2\sum\limits%
_{n=1}^{\infty}q^{n^{2}}\cos\left(2nz\right),\\
\theta_{4}\left(z,q\right)&=1+2\sum\limits%
_{n=1}^{\infty}(-1)^{n}q^{n^{2}}\cos\left(2nz\right).
\end{align*}
We wish to prove that
\begin{equation}
\label{main}
\chi(q^{5})=
\frac{\varphi(q)\, f(q)}{4\, \sum_{n\geq 0}(-q)^{\frac{n(n+1)}{2}}\sin \left\{\frac{(2n+1)3\pi}{10}\right\}\, \sum_{n\geq 0}(-q)^{\frac{n(n+1)}{2}}\sin \left\{\frac{(2n+1)\pi}{10}\right\}}.
\end{equation}
\section{Proof of \eqref{main}}
\begin{proof}
We begin with an identity of  Ramanujan's in his lost notebook \cite[p. 27, Theorem 1.6.1 (ii)]{Bruce2}
$$
\varphi^{2}(q)-5\, \varphi^{2}(q^{5})=-4\, f^{2}(-q^{2})\, \frac{\chi(q^{5})}{\chi(q)}.
$$
The factorizations of the above are given by \cite[p. 29, Entry 1.7.2 (i)]{Bruce2}
\begin{equation}
\label{f1}
\varphi(q)+\sqrt{5}\, \varphi(q^{5})=\frac{(1+\sqrt{5})\, f(-q^{2})}{\prod_{n \text{ odd}}(1+\alpha q^{n}+q^{2n})\, \prod_{n \text{ even}}(1-\beta q^{n}+q^{2n})},
\end{equation}
and \cite[p. 30, Entry 1.7.2 (ii)]{Bruce2}
\begin{equation}
\label{f2}
\varphi(q)-\sqrt{5}\, \varphi(q^{5})=\frac{(1-\sqrt{5})\, f(-q^{2})}{\prod_{n \text{ odd}}(1+\beta q^{n}+q^{2n})\, \prod_{n \text{ even}}(1-\alpha q^{n}+q^{2n})},
\end{equation}
where 
$$
\alpha = \frac{1-\sqrt{5}}{2}, \qquad
\beta=\frac{1+\sqrt{5}}{2}.
$$
Multiplying \eqref{f1} and \eqref{f2} and using the following Fourier expansions \cite[p. 469]{Watson} for $\theta_{1}(z,q)$ and $\theta_{3}(z,q)$:
$$
\theta_{1}\left(z,q\right)=2q^{1/4}\sin z\prod\limits_{n=1}^{\infty}{\left(1-q%
^{2n}\right)}{\left(1-2q^{2n}\cos\left(2z\right)+q^{4n}\right)},
$$
$$
\theta_{3}\left(z,q\right)=\prod\limits_{n=1}^{\infty}\left(1-q^{2n}\right)%
\left(1+2q^{2n-1}\cos\left(2z\right)+q^{4n-2}\right)
$$
we have
$$
\frac{\chi(q^{5})}{\chi(q)}= \frac{q^{1/2}\, f^{4}(-q^{2})}{\theta_{1}(\pi/10,q)\, \theta_{3}(\pi/10,q)\,  \theta_{1}(3\, \pi/10,q)\, \theta_{3}(3\, \pi/10,q)} .
$$
Using the identity \cite{Walker}
$$
\frac{\theta_{1}\left(z,q^{2}\right)\theta_{3}\left(z,q^{2}\right)}{\theta_{1}%
\left(z,iq\right)}=i^{-1/4}\sqrt{\frac{\theta_{2}\left(0,q^%
{2}\right)\theta_{4}\left(0,q^{2}\right)}{2}}
$$
we now have
\begin{align*}
    \frac{\chi(q^{5})}{\chi(q)} &= \frac{2\, q^{1/2}\,i^{1/2} f^{4}(-q^{2})}{\theta_{1}(\pi/10,i\, \sqrt{q})\,\theta_{1}(3\,\pi/10,i\, \sqrt{q})\, \theta_{2}\left(0,q\right)\theta_{4}\left(0,q\right) }\\
    &= \frac{ q^{1/4}\, i^{1/2}\, f^{4}(-q^{2})}{ \psi(q^{2})\, \varphi(-q)\, \theta_{1}(\pi/10,i\, \sqrt{q})\,\theta_{1}(3\,\pi/10,i\, \sqrt{q})}\\
    &= \frac{ q^{1/4}\, i^{1/2}\, f(-q^{2})\, \varphi^{2}(-q^{2})}{ \varphi(-q)\, \theta_{1}(\pi/10,i\, \sqrt{q})\,\theta_{1}(3\,\pi/10,i\, \sqrt{q})}\\
    &= \frac{ q^{1/4}\, i^{1/2}\, f(-q^{2})\, \varphi(q)}{  \theta_{1}(\pi/10,i\, \sqrt{q})\,\theta_{1}(3\,\pi/10,i\, \sqrt{q})}
    \\
    &= \frac{\sum_{n=-\infty}^{\infty}q^{n^{2}}\, \sum_{n=-\infty}^{\infty}(-1)^{n}\, q^{n(3n-1)}}{4\, \sum_{n\geq 0}(-q)^{\frac{n(n+1)}{2}}\sin \left\{\frac{(2n+1)3\pi}{10}\right\}\, \sum_{n\geq 0}(-q)^{\frac{n(n+1)}{2}}\sin \left\{\frac{(2n+1)\pi}{10}\right\}},
\end{align*}
where we have used \cite[p. 15, (1.3.34)]{Bruce1}
$$
f^{3}(-q^{2})=\psi(q^{2})\, \varphi^{2}(-q^{2})
$$
and \cite[p. 15, (1.3.32)]{Bruce1}
$$
\varphi^{2}(-q^{2})=\varphi(q)\, \varphi(-q).
$$
Now we finally use \cite[p.15, (1.3.31)]{Bruce2}
$$
f(q)=\chi(q)\, f(-q^{2})
$$
to conclude our main assertion.
\end{proof}

\end{document}